\documentclass{elsart}
\usepackage{ifpdf}
\usepackage{graphicx,amssymb,lineno}
\usepackage{amssymb}
\usepackage{mathrsfs}
\usepackage{amsmath}
\usepackage{cite}

\newtheorem{theorem}{Theorem}

\newtheorem{remark}{Remark}
\newtheorem{corollary}{Corollary}

\newtheorem{lemma}{Lemma}
\newtheorem{proposition}{Proposition}

\numberwithin{equation}{section} \numberwithin{theorem}{section}
\numberwithin{lemma}{section} \numberwithin{corollary}{section}
\numberwithin{definition}{section}
\numberwithin{proposition}{section} \numberwithin{remark}{section}
\numberwithin{example}{section}
 \ifpdf
\usepackage[%
  pdftitle={Instructions for use of the document class
    elsart},%
  pdfauthor={Simon Pepping},%
  pdfsubject={The preprint document class elsart},%
  pdfkeywords={instructions for use, elsart, document class},%
  pdfstartview=FitH,%
  bookmarks=true,%
  bookmarksopen=true,%
  breaklinks=true,%
  colorlinks=true,%
  linkcolor=blue,anchorcolor=blue,%
  citecolor=blue,filecolor=blue,%
  menucolor=blue,pagecolor=blue,%
  urlcolor=blue]{hyperref}
\else
\usepackage[%
  breaklinks=true,%
  colorlinks=true,%
  linkcolor=blue,anchorcolor=blue,%
  citecolor=blue,filecolor=blue,%
  menucolor=blue,pagecolor=blue,%
  urlcolor=blue]{hyperref}
\fi

\makeatletter
\def\elsartstyle{%
    \def\normalsize{\@setfontsize\normalsize\@xiipt{14.5}}
    \def\small{\@setfontsize\small\@xipt{13.6}}
    \let\footnotesize=\small
    \def\large{\@setfontsize\large\@xivpt{18}}
    \def\Large{\@setfontsize\Large\@xviipt{22}}
    \skip\@mpfootins = 18\p@ \@plus 2\p@
    \normalsize
} \@ifundefined{square}{}{} \makeatother

\begin{document}

\begin{frontmatter}
\title{Blow-up rates  \\for the general curve shortening flow}
\author{RongLi Huang}
\ead{huangronglijane@yahoo.cn}
\footnote{ This work is supported by the National Natural Science
Foundation of China.}
\address{Institute of Mathematics, Fudan University,
\\ Shanghai, 200433, People's Republic of China}

\author{JuanJuan Chen}
\ead{janegirlhappy@yahoo.com.cn}
\address{School of Mathematics and Computational Science, \\Guilin University of Electronic Technology,
\\ Guilin, 541004, People's Republic of China}

\date{}
\begin{abstract}
 \quad The blow-up rates of derivatives of  the curvature function will be presented when the closed curves contract
 to a point in finite time under the general curve shortening flow. In particular,
 this generalizes a theorem of M.E. Gage and R.S. Hamilton about mean curvature  flow in $\mathbb{R}^{2}$.
\end{abstract}

\end{frontmatter}

\section{Introduction}
\quad The curve shortening flow has been studied
extensively in the last thirty years (cf. \cite{B1},\cite{CZ}). As for the newest development  both on  expansion
and contraction of convex closed curves in
$\mathbb{R}^{2},$ see \cite{B2}, \cite{NT}, \cite{LT}, \cite{LPT}, etc.

\quad It is the asymptotic behavior
of the general curve shortening flow that we will mostly be concerned with. Let $\mathbb{S}^{1}$ be an unit circle in the plane, and define
\begin{equation*}
 \gamma_{0}:\,\,\, \mathbb{S}^{1}\rightarrow \mathbb{R}^{2} ,
\end{equation*}
as the closed convex curve in the plane. We look for a  family of closed curves
\begin{equation*}
 \gamma(u,t):\,\,\, \mathbb{S}^{1}\times [0,T)\rightarrow \mathbb{R}^{2} ,
\end{equation*}
which satisfies
\begin{equation}\label{1.1}
\left\{ \begin{aligned} \frac{\partial \gamma}{\partial t}(u,
t)&=G(k)kN , \quad
&u\in \mathbb{S}^{1},\quad &t\in[0, T),\\
\gamma(u,0)&=\gamma_{0}(u), \quad &u\in \mathbb{S}^{1}, \quad &t=0,
\end{aligned} \right.
\end{equation}
where $G$ is a  positive, non-decreasing smooth function on
$(0,\infty)$, $k(\cdot ,t)$ is the inward curvature of the
plane curve $\gamma(\cdot ,t)$ and $N(\cdot ,t)$ is the  unit inward normal
vector.

\quad Assume that $A(t)$ is  the area
of a bounded domain enclosed by the curve $\gamma(\cdot ,t)$,
$L(t)$ is  the length of $\gamma(\cdot, t)$,  $r_{out}(t)$ and $r_{in}(t)$
are respectively the radii of the largest circumscribed circle and the smallest circumscribed circle
of $\gamma(\cdot ,t)$. Define
$$ k_{\max}(t)=\max\{k(u,t)\mid u\in \mathbb{S}^{1}\}, $$
$$ k_{\min}(t)=\min\{k(u,t)\mid u\in \mathbb{S}^{1}\}. $$
\quad The following existence theorem of (\ref{1.1}) belongs to Ben.Andrews (cf. Theorem $\Pi$4.1, Proposition $\Pi$ 4.4 in \cite{B1}).
\begin{proposition}\label{p1.1}
Let $\gamma_{0}$ be a closed strictly convex curve. Then the unique
classical  solution $\gamma(\cdot,t)$ of (\ref{1.1}) exists only at finite time interval $[0,\omega)$, and
 the solution $\gamma(\cdot,t)$ converges to a point $\vartheta$ as
$t\rightarrow \omega$ and $A(t)$, $k_{\max}(t)$  satisfy the
following properties:
\begin{equation*}\label{1.2}
\forall t\in [0,\omega),\, A(t)>0, \,k_{\max}(t)<+\infty,
\end{equation*}
\begin{equation*}\label{1.3}
\lim_{t\rightarrow \omega}A(t)=0, \ \ \ \ \lim_{t\rightarrow \omega}k_{\max}(t)=+\infty.
\end{equation*}
\item As $t\rightarrow \omega$, the normalized curves
$$ \eta(\cdot,t)=\sqrt{\frac{\pi}{A(t)}}\gamma(\cdot,t) $$
converges to the unit circle centered at the point $\vartheta$.

\end{proposition}

\quad Furthermore, if we assume that $G(x)$ is a function on $(0,\infty)$ satisfying

{\bf ($H1$)}\ $G(x) \in C^{3}(0,\infty)$, $G'(x)\geq 0$ and $G(x)>0$ for $x\in (0,\infty)$,

{\bf ($H2$)}\ $G(x)x^{2}$ is convex in $(0,\infty)$ and there is a positive constant $C_{0}$ such that
\begin{equation*}
G'(x)x\leq C_{0}G(x)\qquad \mathrm{for}\,\, \mathrm{sufficiently} \,\, \mathrm{large}\,\, x.
\end{equation*}
Then Rong-Li Huang and Ji-Guang Bao had obtained the following theorem(cf. \cite{HB}).

\begin{proposition}\label{1.2}
Suppose that G(x) satisfies $(H1)$ and $(H2)$.  Let
$\gamma(\cdot,t)$ be the solution for Proposition \ref{p1.1}. Then there
 holds:
\begin{enumerate}

\item
$\displaystyle\lim_{t\rightarrow\omega}\frac{r_{in}(t)}{r_{out}(t)}=1.$

\item
$\displaystyle\lim_{t\rightarrow\omega}\frac{k_{\min}(t)}{k_{\max}(t)}=1.$

\item
 $\displaystyle\lim_{t\rightarrow \omega}\frac{1}{\omega-t}\int_{k(\theta,t)}^{+\infty}\frac{dx}{G(x)x^{3}}=1 $
is uniformly convergent on $\mathbb{S}^{1} $.
\end{enumerate}
Especially, set $G(x)=| x|^{p-1}$ with $p\geq 1$ .  Then
\begin{equation}\label{e1.2}
k(\theta,t)[(p+1)(\omega-t)]^{\frac{1}{p+1}}\,\, \mathrm{converges} \,\, \mathrm{uniformly }\,\,
\mathrm{to}\,\, 1 \,\, \mathrm{as}\,\, t\rightarrow \omega.
\end{equation}
\end{proposition}
\quad Ben.Andrews classifies the limiting shapes for the isotropic curve flows (cf \cite{B2}). However, the blow-up rates
of  derivatives of the curvature function for the general curve shortening flow are not concerned. We now state the main theorem  of this paper.
\begin{theorem}\label{1.3}
Suppose that $G(x)=| x|^{p-1}$ with $p\geq 1$ and $k$ is the curvature function of $\gamma(\cdot,t)$ according to the flow (\ref{1.1}).
$\gamma_{0}$ is a smooth and closed convex curve.
If  $\gamma(\cdot,t)$ converges to a point $\vartheta$ as
$t\rightarrow \omega$.
Then for each $l\in \{1,2,\cdots.\}$, there exists some positive constant $C(l,p)$
depending only on $l$, $p$, such that
\begin{equation}\label{1.5}
\parallel\frac{\partial^{l}k}{\partial\theta^{l}}\parallel_{L^{\infty}(\mathbb{S}^{1})}\leq C(l,p)(\omega-t)^{2\alpha
\frac{p}{p+1}-\frac{1}{p+1}},
\end{equation}
where $\alpha$ is any constant satisfying $0<\alpha<1$ .
\end{theorem}

\begin{remark}\label{1.4}
By (\ref{1.5}) we  generalize  Corollary 5.7.2 in \cite{GH}.
\end{remark}

\quad This paper is organized as follows: In the next section we obtain an equivalent proposition of Theorem \ref{1.3}
by means of transferring the flow (\ref{1.1}) into an Cauchy PDE's problem. Finally in the third section,
it proves Proposition \ref{p2.1} by making use of
Gage-Hamilton's method.

\section{Preliminaries}
\quad By Lemma 2.4 in \cite{HB}, the general curve shortening problem
(\ref{1.1}) with initial convex curve $\gamma_{0}(\theta)$ is equivalent to the cauchy problem
\begin{equation}\label{e2.1}
\left\{ \begin{aligned} \frac{\partial k}{\partial
t}&=k^{2}\left(\frac{\partial^{2}}{\partial\theta^{2}}(G(k)k)+G(k)k\right),\quad
&\theta \in \mathbb{S}^{1},\quad &t\in(0, \omega),\\
k&=k_{0}(\theta),
&\theta\in\mathbb{S}^{1},\quad &t=0.
\end{aligned} \right.
\end{equation}
where $0<\alpha<1$, $k\in
C^{2+\alpha,1+\frac{\alpha}{2}}(\mathbb{S}^{1}\times (0,\omega))$,
$k_{0}(\theta)$ denotes the curvature function of the curve
$\gamma_{0}(\theta)$.

\quad Assume that  $G(x)=| x|^{p-1}$ with $p\geq 1$. Set
\begin{equation*}
\tau=-\frac{1}{p+1}\ln\frac{\omega-t}{\omega},\quad\tilde{k}(\theta,\tau)=k^{p}(\theta,t)[(p+1)(\omega-t)]^{\frac{p}{P+1}}.
\end{equation*}
It follows from (\ref{e2.1}) that $\tilde{k}(\theta,\tau)$ satisfies
\begin{equation}\label{e2.2}
\left\{ \begin{aligned} \frac{\partial \tilde{k}}{\partial
\tau}&=p\tilde{k}\tilde{k}^{\frac{1}{p}}\frac{\partial^{2}\tilde{k}}{\partial\theta^{2}}+p\tilde{k}^{2}\tilde{k}^{\frac{1}{p}}-p\tilde{k},\quad
&\theta \in \mathbb{S}^{1},\quad &\tau\in(0, +\infty),\\
\tilde{k}&=\tilde{k}_{0}(\theta),
&\theta\in\mathbb{S}^{1},\quad &\tau=0,
\end{aligned} \right.
\end{equation}
where $\tilde{k}_{0}(\theta)=k_{0}^{p}(\theta)((p+1)\omega)^{\frac{p}{P+1}}$ is a positive smooth function. By (\ref{e1.2}), there holds
\begin{equation}\label{e2.3}
\lim_{\tau\rightarrow+\infty}\tilde{k}(\theta,\tau)=1 \,\,\,\mathrm{uniformly} \,\,\mathrm{on} \,\,\mathbb{S}^{1}.
\end{equation}
\quad Define
\begin{equation*}
\|\tilde{k}^{(l)}\|_{q}=\biggl[\int_{\mathbb{S}^{1}}\biggl|\frac{\partial^{l}\tilde{k}}{\partial\theta^{l}}\biggl|^{q}\biggr]^{\frac{1}{q}}
\qquad \mathrm{for}\quad l\in \{1,2,\cdots.\} \quad\mathrm{and}\quad q\geq 1.
\end{equation*}
Here,  $\tilde{k}'$, $\tilde{k}'',\cdots, \tilde{k}^{(l)},\cdots,$ denote   partial differentiation by $\theta$.
 Then Theorem \ref{1.3} is equivalent to the following proposition.
\begin{proposition}\label{p2.1}
If $\tilde{k}(\theta,\tau)$ is a smooth solution of  problem (\ref{e2.2}) with $p\geq 1$ and  satisfies (\ref{e2.3}). Then for each $l\in \{1,2,\cdots.\}$ there exists some constant $C(l,p)$ depending only on $l$, $p$  and the initial curve $\gamma_{0}$  such that
\begin{equation}\label{e2.4}
\|\tilde{k}^{(l)}\|_{\infty}\leq C(l,p)\exp(-2\alpha p\tau),\quad \forall \tau>0,
\end{equation}
where $\alpha$ is any constant satisfying $0<\alpha<1$.
\end{proposition}

\quad Obviously our task in the following can be changed over to prove Proposition \ref{p2.1}.
The forging listing of three facts and two lemmas can be
used repeatedly (cf. \cite{GH}).

$(\mathbb{I})$. Young inequality. For all positive $\epsilon$, $p$, $q$,
 $\displaystyle ab\leq \epsilon^{p} \frac{a^{p}}{p}+\frac{b^{q}}{\epsilon^{q}q}$ where $\displaystyle\frac{1}{p}+\frac{1}{q}=1$.

$(\mathbb{II})$. Wirtinger inequality. If $\displaystyle\int_{\mathbb{S}^{1}}f=0$,
then $\displaystyle\int_{\mathbb{S}^{1}}f^{2}\leq\int_{\mathbb{S}^{1}}f'^{2}$.

$(\mathbb{III})$. Sobolev inequality. If $\parallel f\parallel_{2}\leq C$ and $\parallel f'\parallel_{2}\leq C$,
then  $\displaystyle\parallel f\parallel_{\infty}\leq\biggl(\frac{1}{\sqrt{2\pi}}+\sqrt{2\pi}\biggr) C$.
\begin{lemma}[Gage-Hamilton] \label{l2.2}
Let $f: \mathbf{R}^{+}\rightarrow\mathbf{R}^{+}$ satisfy
\begin{equation*}
\displaystyle\frac{df}{d\tau}\leq Cf^{1-\frac{1}{q}}-2C(q)f\quad\, \mathrm{where} \quad\,q\geq 1.
\end{equation*}
Then
\begin{equation*}
 \displaystyle f^\frac{1}{q}(\tau)\leq \biggl(\frac{C}{2q}+D\exp(-2\tau)\biggr)\leq\tilde{C}(q).
 \end{equation*}
\end{lemma}
\begin{lemma}[Gage-Hamilton]\label{l2.3}
If $\displaystyle\frac{df}{d\tau}\leq-\alpha f+C\exp(-\beta\tau)$, then
\begin{equation*}
f(\tau)\leq D\exp(-\alpha\tau)+\frac{C}{\alpha-\beta}\exp(-\beta\tau) \quad \mathrm{if}\quad \alpha\neq\beta,\quad \mathrm{or}
\end{equation*}
\begin{equation*}
\,\,\,f(\tau)\leq D\exp(-\alpha\tau)+C\tau\exp(-\alpha\tau) \quad \mathrm{if}\quad \alpha=\beta.\qquad\qquad
\end{equation*}
\end{lemma}

\section{ Decay estimates of the scaling flow}
\quad In this section we establish some formal estimates on the analogy of \cite{GH} in order to obtain (\ref{e2.4}). In addition,
those heuristic deductions bring out the decay estimates of the scaling curvature function.
Hereafter we always suppose that $\tilde{k}$ is the smooth positive solution of  (\ref{e2.2}) and satisfies (\ref{e2.3}).
\begin{lemma}\label{l3.1}
There exists a constant $C_{1}$ depending only on $p$ and the initial curve $\gamma_{0}$ such that there holds
\begin{equation}\label{e3.1}
\parallel\tilde{k}'\parallel_{2}\leq C_{1},\quad \parallel\tilde{k}'\parallel_{4}\leq C_{1}.
\end{equation}
\end{lemma}
{\bf Proof.}
For $n\in \{1,3,5,7,\cdots,\}$, we write
\begin{equation*}
f(\tau)=\int_{\mathbb{S}^{1}}\tilde{k}'^{n+1}.
\end{equation*}
In view of the equation (\ref{l2.2}),  it gives
\begin{equation} \begin{aligned}\label{e3.2}
\frac{f'(\tau)}{(n+1)p}&=\int_{\mathbb{S}^{1}}\tilde{k}'^{n}(\tilde{k}\tilde{k}^{\frac{1}{p}}\tilde{k}''
+\tilde{k}^{2}\tilde{k}^{\frac{1}{p}}-\tilde{k})'\\
&=-\int_{\mathbb{S}^{1}}\tilde{k}'^{n+1}-n\int_{\mathbb{S}^{1}}\tilde{k}'^{n-1}\tilde{k}''^{2}\tilde{k}^{\frac{p+1}{p}}\\
&\quad-n\int_{\mathbb{S}^{1}}\tilde{k}'^{\frac{n-1}{2}}\tilde{k}''\tilde{k}^{\frac{p+1}{2p}}\tilde{k}'^{\frac{n-1}{2}}\tilde{k}^{\frac{3p+1}{2p}}\\
&\leq-f(\tau)-n\int_{\mathbb{S}^{1}}\tilde{k}'^{n-1}\tilde{k}''^{2}\tilde{k}^{\frac{p+1}{p}}\\
&\quad+n\int_{\mathbb{S}^{1}}\tilde{k}'^{n-1}\tilde{k}''^{2}\tilde{k}^{\frac{p+1}{p}}
+\frac{n}{4}\int_{\mathbb{S}^{1}}\tilde{k}'^{n-1}\tilde{k}^{\frac{3p+1}{p}}\\
&=-f(\tau)+\frac{n}{4}\int_{\mathbb{S}^{1}}\tilde{k}'^{n-1}\tilde{k}^{\frac{3p+1}{p}}\\
&\leq -f(\tau)+\epsilon^{\frac{n+1}{n-1}}\frac{n}{4}\frac{n-1}{n+1}\int_{\mathbb{S}^{1}}\tilde{k}'^{n+1}\\
&\quad+\epsilon^{-\frac{n+1}{2}}\frac{n}{4}\frac{2}{n+1}\int_{\mathbb{S}^{1}}\tilde{k}^{\frac{n+1}{2}\frac{3p+1}{p}}
\end{aligned}
\end{equation}
in terms of calculating the derivative of $f(\tau)$ and integrating by parts.

We take $\epsilon$ such that
\begin{equation*}
\epsilon^{\frac{n+1}{n-1}}\frac{n}{4}\frac{n-1}{n+1}=\frac{1}{2}
\end{equation*}
and use (\ref{e2.3}) then (\ref{e3.2}) implies
\begin{equation*}
f'(\tau)\leq-\frac{(n+1)p}{2}f(\tau)+C(n,p).
\end{equation*}
By Lemma \ref{l2.2}  we get the upper bound
\begin{equation*}
f(\tau)\leq C(n,p),
\end{equation*}
and this yields (\ref{e3.1}) for $n=1, n=3$.
\qed
\begin{lemma}\label{l3.2}
There exists some constant $C_{2}$ depending only on $p$ and the initial curve $\gamma_{0}$ such that
\begin{equation}\label{e3.3}
\parallel\tilde{k}''\parallel_{2}\leq C_{2}.
\end{equation}
\end{lemma}
{\bf Proof.}
Define
\begin{equation*}
g(\tau)=\int_{\mathbb{S}^{1}}\tilde{k}''^{2}.
\end{equation*}
Owing to (\ref{e2.3}) and by analogy with (\ref{e3.2})  we have
\begin{equation} \begin{aligned}\label{e3.4}
\frac{g'(\tau)}{2p}&=\int_{\mathbb{S}^{1}}\tilde{k}''(\tilde{k}\tilde{k}^{\frac{1}{p}}\tilde{k}''
+\tilde{k}^{2}\tilde{k}^{\frac{1}{p}}-\tilde{k})''\\
&=-\int_{\mathbb{S}^{1}}\tilde{k}''^{2}-\int_{\mathbb{S}^{1}}\tilde{k}'''^{2}\tilde{k}^{\frac{p+1}{p}}\\
&\quad-\frac{p+1}{p}\int_{\mathbb{S}^{1}}\tilde{k}'''\tilde{k}^{\frac{p+1}{2p}}\tilde{k}''\tilde{k}'\tilde{k}^{\frac{1-p}{2p}}
-\frac{2p+1}{p}\int_{\mathbb{S}^{1}}\tilde{k}'''\tilde{k}^{\frac{p+1}{2p}}\tilde{k}'\tilde{k}^{\frac{p+1}{2p}}
\\
&\leq-g(\tau)-\int_{\mathbb{S}^{1}}\tilde{k}'''^{2}\tilde{k}^{\frac{p+1}{p}}\\
&\quad+\frac{1}{2}\epsilon^{2}_{1}\frac{p+1}{p}
\int_{\mathbb{S}^{1}}\tilde{k}'''^{2}\tilde{k}^{\frac{p+1}{p}}+\frac{1}{2}\epsilon^{-2}_{1}\frac{p+1}{p}
\int_{\mathbb{S}^{1}}\tilde{k}''^{2}\tilde{k}'^{2}\tilde{k}^{\frac{1-p}{p}}\\
&\quad+\frac{1}{2}\epsilon^{2}_{2}\frac{2p+1}{p}\int_{\mathbb{S}^{1}}\tilde{k}'''^{2}\tilde{k}^{\frac{p+1}{p}}
+\frac{1}{2}\epsilon^{-2}_{2}\frac{2p+1}{p}\int_{\mathbb{S}^{1}}\tilde{k}'^{2}\tilde{k}^{\frac{p+1}{p}}.
\end{aligned}
\end{equation}
Choose $\epsilon_{1}, \epsilon_{2}$ such that
\begin{equation*}
\frac{1}{2}\epsilon^{2}_{1}\frac{p+1}{p}=\frac{1}{2},\quad \frac{1}{2}\epsilon^{2}_{2}\frac{2p+1}{p}=\frac{1}{2}.
\end{equation*}
Consequently, putting (\ref{e2.4}), (\ref{e3.4}) and Lemma \ref{l3.1} together, it shows
\begin{equation}\label{e3.5}
g'(\tau)\leq-2pg(\tau)+C_{1}(p)\int_{\mathbb{S}^{1}}\tilde{k}''^{2}\tilde{k}'^{2}+C_{2}(p)
\end{equation}
From  (\ref{e3.2}) we see that for $n=3$ there holds
\begin{equation}\label{e3.6}
\begin{aligned}
3\int_{\mathbb{S}^{1}}\tilde{k}'^{2}\tilde{k}''^{2}\tilde{k}^{\frac{p+1}{p}}&=
-\frac{1}{4p}\frac{\partial}{\partial\tau}\int_{\mathbb{S}^{1}}\tilde{k}'^{4}
-\int_{\mathbb{S}^{1}}\tilde{k}'^{4}
-3\int_{\mathbb{S}^{1}}\tilde{k}'^{2}\tilde{k}''\tilde{k}^{\frac{2p+1}{2p}}\\
&\leq -\frac{1}{4p}\frac{\partial}{\partial\tau}\int_{\mathbb{S}^{1}}\tilde{k}'^{4}
+\frac{3}{2}\epsilon_{3}^{2}\int_{\mathbb{S}^{1}}\tilde{k}''^{2}+
\frac{3}{2}\epsilon_{3}^{-2}\int_{\mathbb{S}^{1}}\tilde{k}'^{4}\tilde{k}^{\frac{2p+1}{p}}.
\end{aligned}
\end{equation}
By Lemma 2.6 in \cite{HB} and (\ref{e2.3}), $\tilde{k}$ has positive lower bound and upper bound. Thus combining  (\ref{e3.5}) with (\ref{e3.6}) we deduce that
\begin{equation*}\begin{aligned}
g'(\tau)&\leq-2pg(\tau)-C_{3}(p)\frac{\partial}{\partial\tau}\int_{\mathbb{S}^{1}}\tilde{k}'^{4}
+\epsilon_{3}^{2}C_{4}(p)\int_{\mathbb{S}^{1}}\tilde{k}''^{2}\\
&\quad +\epsilon_{3}^{-2}C_{5}(p)\int_{\mathbb{S}^{1}}\tilde{k}'^{4}\tilde{k}^{\frac{2p+1}{p}}+C_{6}(p).
\end{aligned}
\end{equation*}
Taking $\displaystyle\epsilon_{3}=\sqrt{\frac{p}{C_{4}(p)}}$ and using (\ref{e3.1}) we obtain
\begin{equation}\label{e3.7}
g'(\tau)\leq-pg(\tau)-C_{3}(p)\frac{\partial}{\partial\tau}\int_{\mathbb{S}^{1}}\tilde{k}'^{4}
+C_{7}(p).
\end{equation}
Then the inequality (\ref{e3.7}) implies
\begin{equation*}
\frac{\partial}{\partial\tau}[e^{p\tau}g(\tau)]\leq-C_{3}(p)e^{p\tau}\frac{\partial}{\partial\tau}\int_{\mathbb{S}^{1}}\tilde{k}'^{4}
+C_{7}(p)e^{p\tau}.
\end{equation*}
Integrating from $0$ to $A$ on both sides and using (\ref{e3.1}) we find
\begin{equation*}
e^{pA}g(A)\leq C_{8}(p)e^{pA}+C_{9}(p, \gamma_{0}).
\end{equation*}
So that this establish (\ref{e3.3}).
\qed

\quad Similarly, using the methods from the proof of Lemma 5.7.8 in \cite{GH}, we have
\begin{corollary}\label{c3.3}
$\parallel\tilde{k}'\parallel_{\infty}$ converges to zero as $\tau\rightarrow+\infty$.
\end{corollary}
\quad
It is convenient to introduce the following proposition  (cf. Lemma 5.7.9 in \cite{GH}) which will provide us with
a key means to obtain the  estimates  (\ref{e2.4}).
\begin{proposition}[Gage-Hamilton]\label{p3.1}
For any $0<\alpha<1$ we can choose $A$ so that for $\tau\geq A$
\begin{equation*}
4\alpha\int_{\mathbb{S}^{1}}\tilde{k}'^{2}\leq\int_{\mathbb{S}^{1}}\tilde{k}''^{2}.
\end{equation*}
\end{proposition}
\quad  We can next obtain the exponential decay estimates which are similar to Lemma 5.7.10 in \cite{GH}.
\begin{lemma}
For any $\alpha$, $0<\alpha<1$, there is a constant $C_{3}$ depending only on $p$ and the initial curve $\gamma_{0}$ such that
\begin{equation}\label{e3.8}
\parallel\tilde{k}'\parallel_{2}\leq C_{3}\exp(-2\alpha p\tau).
\end{equation}
\end{lemma}
{\bf Proof.}
Without loss we may assume $p>1$ and  $\displaystyle\frac{1}{2}(1+\frac{1}{p})<\alpha<1$. As same as (\ref{e3.2}) we calculate that
\begin{equation} \begin{aligned}\label{e3.9}
\frac{\frac{\partial}{\partial\tau}\int_{\mathbb{S}^{1}}\tilde{k}'^{2}}{2p}&
=\int_{\mathbb{S}^{1}}\tilde{k}'(\tilde{k}\tilde{k}^{\frac{1}{p}}\tilde{k}''
+\tilde{k}^{2}\tilde{k}^{\frac{1}{p}}-\tilde{k})'\\
&=-\int_{\mathbb{S}^{1}}\tilde{k}'^{2}-\int_{\mathbb{S}^{1}}\tilde{k}''^{2}\tilde{k}^{\frac{p+1}{p}}
+\frac{2p+1}{p}\int_{\mathbb{S}^{1}}\tilde{k}'^{2}\tilde{k}^{\frac{p+1}{p}}.
\end{aligned}
\end{equation}
It follows from (\ref{e2.3}) and Proposition \ref{p3.1} that  there is a positive constant $A$  so large that
for $\tau\geq A$ we have
\begin{equation}\label{e3.10}
1-\epsilon\leq\tilde{k}^{\frac{p+1}{p}}\leq 1+\epsilon ,\quad
4\alpha\int_{\mathbb{S}^{1}}\tilde{k}'^{2}\leq\int_{\mathbb{S}^{1}}\tilde{k}''^{2}
\end{equation}
where
\begin{equation*}
\epsilon=\frac{2\alpha-1-\frac{1}{p}}{2+\frac{1}{p}+4\alpha}.
\end{equation*}
Combining  (\ref{e3.9}) with (\ref{e3.10}) for $\tau\geq A$ we obtain
\begin{equation*} \begin{aligned}
\frac{\frac{\partial}{\partial\tau}\int_{\mathbb{S}^{1}}\tilde{k}'^{2}}{2p}&
\leq-\int_{\mathbb{S}^{1}}\tilde{k}'^{2}-(1-\epsilon)4\alpha\int_{\mathbb{S}^{1}}\tilde{k}'^{2}
+\frac{2p+1}{p}(1+\epsilon)\int_{\mathbb{S}^{1}}\tilde{k}'^{2}\\
&=-2\alpha\int_{\mathbb{S}^{1}}\tilde{k}'^{2}.
\end{aligned}
\end{equation*}
By Gronwall's inequality the assertion follows easily.
\qed

\quad Further, we have the  conclusions below.
\begin{lemma}\label{l3.6}
For any $\alpha$, $0<\alpha<1$, there is a constant $C_{4}$ depending only on $p$ and the initial curve $\gamma_{0}$ such that
\begin{equation}\label{e3.11}
\parallel\tilde{k}''\parallel_{2}\leq C_{4}\exp(-2\alpha p\tau).
\end{equation}
\end{lemma}
{\bf Proof.}
Similar to the calculation of (\ref{e3.2}), we arrive at
\begin{equation} \begin{aligned}\label{e3.12}
\frac{\frac{\partial}{\partial\tau}\int_{\mathbb{S}^{1}}\tilde{k}''^{2}}{2p}&
=\int_{\mathbb{S}^{1}}\tilde{k}''(\tilde{k}\tilde{k}^{\frac{1}{p}}\tilde{k}''
+\tilde{k}^{2}\tilde{k}^{\frac{1}{p}}-\tilde{k})''\\
&=-\int_{\mathbb{S}^{1}}\tilde{k}''^{2}-\int_{\mathbb{S}^{1}}\tilde{k}'''^{2}\tilde{k}^{\frac{p+1}{p}}\\
&\quad-\frac{p+1}{p}\int_{\mathbb{S}^{1}}\tilde{k}'''\tilde{k}''\tilde{k}'\tilde{k}^{\frac{1}{p}}-(2+\frac{1}{p})
\int_{\mathbb{S}^{1}}\tilde{k}'''\tilde{k}'\tilde{k}^{\frac{p+1}{p}}
\end{aligned}
\end{equation}
Using (\ref{e2.3}),  (\ref{e3.12}), Corollary \ref{c3.3} and letting $\tau$ so large  that we have the estimate
\begin{equation*} \begin{aligned}
\frac{\frac{\partial}{\partial\tau}\int_{\mathbb{S}^{1}}\tilde{k}''^{2}}{2p}
&\leq-\int_{\mathbb{S}^{1}}\tilde{k}''^{2}-(1-\frac{\epsilon}{2})\int_{\mathbb{S}^{1}}\tilde{k}'''^{2}
+\frac{\epsilon}{4}\int_{\mathbb{S}^{1}}\tilde{k}'''^{2}+\epsilon\int_{\mathbb{S}^{1}}\tilde{k}''^{2}\\
&\quad+\frac{\epsilon}{4}\int_{\mathbb{S}^{1}}\tilde{k}'''^{2}+ \frac{C_{3}}{\epsilon}\exp(-4\alpha p\tau)\\
&\leq-(1-\epsilon)\int_{\mathbb{S}^{1}}\tilde{k}''^{2}-(1-\epsilon)\int_{\mathbb{S}^{1}}\tilde{k}'''^{2}
+ \frac{C_{3}}{\epsilon}\exp(-4\alpha p\tau)\\
&\leq -(2-2\epsilon)\int_{\mathbb{S}^{1}}\tilde{k}''^{2}
+ \frac{C_{3}}{\epsilon}\exp(-4\alpha p\tau)
\end{aligned}
\end{equation*}
 provided any positive constants $\epsilon<1$. Here we use Young inequality and Wirtinger inequality.
 Taking $\epsilon=1-\alpha$ we obtain
\begin{equation*} \begin{aligned}
\frac{\frac{\partial}{\partial\tau}\int_{\mathbb{S}^{1}}\tilde{k}''^{2}}{2p}
&\leq-2\alpha\int_{\mathbb{S}^{1}}\tilde{k}''^{2}+ \frac{C_{3}}{1-\alpha}\exp(-4\alpha p\tau)
\end{aligned}
\end{equation*}

Furthermore, based on Lemma \ref{l2.3} ,  the assertion follows.
\qed

\quad Applying Sobolev inequality to (\ref{e3.8}) and (\ref{e3.11}) implies the following conclusions.
\begin{corollary}\label{c3.7}
For any $\alpha$, $0<\alpha<1$, there is a constant $C_{5}$ depending only on $p$ and
the initial curve $\gamma_{0}$ such that
\begin{equation*}
\parallel\tilde{k}'\parallel_{\infty}\leq C_{5}\exp(-2\alpha p\tau).
\end{equation*}
\end{corollary}
{ \bf Proof of  Proposition 2.1 .}

\quad We will give the complete proof of this proposition by Mathematical Induction.
First we get
\begin{equation*}
\|\tilde{k}^{(L)}\|_{\infty}\leq C(L,p)\exp(-2\alpha p\tau),\quad \forall \tau>0,
\end{equation*}
 and
\begin{equation*}
\parallel\tilde{k}^{(L+1)}\parallel_{2}\leq C(L,p)\exp(-2\alpha p\tau), \quad \forall \tau>0,
\end{equation*}
for $L=1$ by Corollary \ref{c3.7} and Lemma \ref{l3.6}. Next we assume that
\begin{equation}\label{e3.13}
\|\tilde{k}^{(l)}\|_{\infty}\leq C(l,p)\exp(-2\alpha p\tau),\quad \forall \tau>0,
\end{equation}
and
\begin{equation}\label{e3.14}
\parallel\tilde{k}^{(l+1)}\parallel_{2}\leq C(l,p)\exp(-2\alpha p\tau), \quad \forall \tau>0,
\end{equation}
for $l=1, 2, \cdots, L$ provided $L>1$.  If  we can deduce that
\begin{equation}\label{e3.16}
\|\tilde{k}^{(L+1)}\|_{\infty}\leq C(L+1,p)\exp(-2\alpha p\tau),\quad \forall \tau>0,
\end{equation}
by the assumption of (\ref{e3.13}) and (\ref{e3.14}) then we obtain the desired results.

\quad To prove (\ref{e3.16}) we need only
\begin{equation}\label{e3.17}
\|\tilde{k}^{(L+2)}\|_{2}\leq C(L+1,p)\exp(-2\alpha p\tau),\quad \forall \tau>0,
\end{equation}
by  the Sobolev inequality.

\quad In order to verify (\ref{e3.17}) we show that
\begin{equation} \begin{aligned}\label{e3.18}
\frac{\frac{\partial}{\partial\tau}\int_{\mathbb{S}^{1}}[\tilde{k}^{(L+2)}]^{2}}{2p}&
=\int_{\mathbb{S}^{1}}\tilde{k}^{(L+2)}(\tilde{k}\tilde{k}^{\frac{1}{p}}\tilde{k}''
+\tilde{k}^{2}\tilde{k}^{\frac{1}{p}}-\tilde{k})^{(L+2)}\\
&=-\int_{\mathbb{S}^{1}}[\tilde{k}^{(L+2)}]^{2}-\int_{\mathbb{S}^{1}}[\tilde{k}^{(L+3)}]^{2}\tilde{k}^{\frac{p+1}{p}}\\
&\quad-\sum_{i=1}^{L}\int_{\mathbb{S}^{1}}\tilde{k}^{(L+3)}\tilde{k}^{(L+1-i)}(\tilde{k}^{1+\frac{1}{p}})^{(i)}\\
&\quad-\int_{\mathbb{S}^{1}}\tilde{k}^{(L+3)}(\tilde{k}^{1+\frac{1}{p}})^{(L+1)}
-\int_{\mathbb{S}^{1}}\tilde{k}^{(L+3)}(\tilde{k}^{2+\frac{1}{p}})^{(L+1)}.
\end{aligned}
\end{equation}
Let $\epsilon$ be small constant to be determined and take $\tau$ to be large enough. Then by (\ref{e2.3})  we have
\begin{equation}\label{e3.19}
\int_{\mathbb{S}^{1}}[\tilde{k}^{(L+3)}]^{2}\tilde{k}^{\frac{p+1}{p}}
\geq (1-\epsilon)\int_{\mathbb{S}^{1}}[\tilde{k}^{(L+3)}]^{2}.
\end{equation}
Hence by (\ref{e3.13}) and the Young inequality,
\begin{equation}\begin{aligned}\label{e3.20}
&-\sum_{i=1}^{L}\int_{\mathbb{S}^{1}}\tilde{k}^{(L+3)}\tilde{k}^{(L+1-i)}(\tilde{k}^{1+\frac{1}{p}})^{(i)}\\
&\leq \frac{\epsilon}{2}\int_{\mathbb{S}^{1}}[\tilde{k}^{(L+3)}]^{2}+\frac{C(L,p)}{\epsilon}\exp(-4\alpha p\tau).
\end{aligned}
\end{equation}
By (\ref{e3.14}),  using the Young inequality again we obtain
\begin{equation}\begin{aligned}\label{e3.21}
&-\int_{\mathbb{S}^{1}}\tilde{k}^{(L+3)}(\tilde{k}^{1+\frac{1}{p}})^{(L+1)}\leq  \frac{\epsilon}{4}\int_{\mathbb{S}^{1}}[\tilde{k}^{(L+3)}]^{2}+\frac{C(L,p)}{\epsilon}\exp(-4\alpha p\tau),\\
&-\int_{\mathbb{S}^{1}}\tilde{k}^{(L+3)}(\tilde{k}^{2+\frac{1}{p}})^{(L+1)}\leq  \frac{\epsilon}{4}\int_{\mathbb{S}^{1}}[\tilde{k}^{(L+3)}]^{2}+\frac{C(L,p)}{\epsilon}\exp(-4\alpha p\tau).
\end{aligned}
\end{equation}
Substituting (\ref{e3.19})-(\ref{e3.21}) to (\ref{e3.18}) and using the Wirtinger inequality we get
\begin{equation*} \begin{aligned}
\frac{\frac{\partial}{\partial\tau}\int_{\mathbb{S}^{1}}[\tilde{k}^{(L+2)}]^{2}}{2p}&\leq -\int_{\mathbb{S}^{1}}[\tilde{k}^{(L+2)}]^{2}
-(1-2\epsilon)\int_{\mathbb{S}^{1}}[\tilde{k}^{(L+3)}]^{2}+\frac{C(L,p)}{\epsilon}\exp(-4\alpha p\tau)\\
&\leq -\int_{\mathbb{S}^{1}}[\tilde{k}^{(L+2)}]^{2}
-(1-2\epsilon)\int_{\mathbb{S}^{1}}[\tilde{k}^{(L+2)}]^{2}+\frac{C(L,p)}{\epsilon}\exp(-4\alpha p\tau)\\
&= -(2-2\epsilon)\int_{\mathbb{S}^{1}}[\tilde{k}^{(L+2)}]^{2}
+\frac{C(L,p)}{\epsilon}\exp(-4\alpha p\tau).
\end{aligned}
\end{equation*}
Choose $\epsilon=1-\alpha$ then there holds
\begin{equation}\label{e3.22}
\frac{\frac{\partial}{\partial\tau}\int_{\mathbb{S}^{1}}[\tilde{k}^{(L+2)}]^{2}}{2p}\leq
-2\alpha\int_{\mathbb{S}^{1}}[\tilde{k}^{(L+2)}]^{2}
+\frac{C(L,p)}{1-\alpha}\exp(-4\alpha p\tau).
\end{equation}
Applying Lemma \ref{l2.3} to (\ref{e3.22}) we obtain (\ref{e3.17}). This completes the proof.
\qed

\end{document}